\input ./style/arxiv-general.cfg
\documentclass[MSNbibl,number,citesort,seceqn,dvips]{arxbj}
\makeatletter
   \@ifpackageloaded{graphicx}{}{\usepackage{graphicx}}
\makeatother

\aid{0}
\volume{21}
\issue{3}
\pubyear{2015}
\firstpage{1844}
\lastpage{1854}
\doi{10.3150/14-BEJ628} 

\makeatletter
\newtheorem{them}{Theorem}[section]
\newproclaim{Notation}{Notation}
\newproclaim{Definition}{Definition}

\newcommand{\eqref}[1]{(\ref{#1})}
\renewcommand{\epsilon}{\varepsilon}
\renewcommand{\emptyset}{\varnothing}
\newcommand{\und}{\underline}
\newcommand{\La}{\Lambda}
\newcommand{\eps}{\epsilon}
\newcommand{\ga}{\gamma}
\newcommand{\Om}{\Omega}
\makeatother

\begin{document}
\begin{frontmatter}

\title{Exponential rate of convergence in current~reservoirs}
\runtitle{Spectral gap}

\begin{aug}
\author[1]{\inits{A.}\fnms{Anna}~\snm{De Masi}\corref{}\thanksref{1}\ead[label=e1]{demasi@univaq.it}},
\author[2]{\inits{E.}\fnms{Errico}~\snm{Presutti}\thanksref{2}\ead[label=e2]{errico.presutti@gmail.com}},
\author[3]{\inits{D.}\fnms{Dimitrios}~\snm{Tsagkarogiannis}\thanksref{3}\ead[label=e3]{D.Tsagkarogiannis@sussex.ac.uk}} \and
\author[4]{\inits{M.E.}\fnms{Maria Eulalia}~\snm{Vares}\thanksref{4}\ead[label=e4]{eulalia@im.ufrj.br}}
\address[1]{Dipartimento di Ingegneria e Scienze dell'Informazione e
Matematica, Universit\`a di L'Aquila, Via Vetoio, 1 67100 L'Aquila,
Italy. \printead{e1}}
\address[2]{Gran Sasso Science Institute, Viale Francesco Crispi, 7
67100 L'Aquila, Italy.\\ \printead{e2}}
\address[3]{Department of Mathematics, University
of Sussex, Pevensey 2 Building, Falmer Campus, Brighton BN1 9QH, UK. \printead{e3}}
\address[4]{Instituto de Matem\'{a}tica, Universidade Federal do Rio
de Janeiro, Av. Athos da S. Ramos 149, 21941-909, Rio de Janeiro, RJ,
Brazil. \printead{e4}}
\end{aug}

\received{\smonth{3} \syear{2013}}
\revised{\smonth{12} \syear{2013}}

%
\begin{abstract}
In this paper, we consider a family of interacting particle systems on
$[-N, N]$ that arises as a natural model for
current reservoirs and Fick's law. We study the exponential rate of
convergence to the stationary measure, which we prove to be of the
order $N^{-2}$.
\end{abstract}

%
\end{frontmatter}

\section{Introduction}
\label{sec:z1}

In this paper, we study a family of interacting particle systems whose
state space is $\{0,1\}^{[-N,N]}$. For each $N$, the
dynamics is a Markov process with generator $L=L_0+L_b$, $L_0$ the
generator of the stirring process (see \eqref{z2.1} below),
$L_b$ the generator of a birth-death process whose events are localized
in a neighborhood of the end-points; see \eqref{z3.1}.

In particular, we focus on the case when around $N$ there are only
births while around $-N$ there are only
deaths. The system is then ``unbalanced'' and in the stationary measure
$\mu^{\mathrm{st}}_N$ there is a non-zero steady current
of particles flowing from right to left. This system is designed to
model the Fick's law which relates the current
to the density gradient.

In statistical mechanics, non-equilibrium is not as well understood as
equilibrium, hence the interest,
from a physical viewpoint, to look at systems which are stationary yet
in non-equilibrium: in our case,
the stationary process is in fact non-reversible and the stationary
measure $\mu^{\mathrm{st}}_N$ not Gibbsian.

There is a huge literature on stationary non-equilibrium measures, in
particular,
on their large deviations, as they are related to ``out of equilibrium
thermodynamics''
(see, for instance, \mbox{\cite{bertini,bodineau,bdl,derrida}}). Our goal
is to study the exponential rate at which the dynamics converges to the
stationary measure, and how it depends on the system size. Spectral
gaps have been well studied in the reversible or Gibbsian set-up, both
for stirring and for more general interacting particle systems (see,
for instance, \cite{lu}). The techniques\vadjust{\goodbreak} used in those situations,
however, do not seem to apply to our non-equilibrium model. We shall
rather rely on stochastic inequalities and coupling methods, thus reducing
the problem to that of bounding the extinction time of the set of
discrepancies between two coupled evolutions.
The case of a single discrepancy can be regarded as an environment
dependent random walk with death rate
which also depends on the environment. Its extinction time has been
studied in \cite{dptvrw} and, as we shall see here,
is closely related to the exponential rate of convergence in our model.

The main part of this paper refers to the case of ``current reservoirs''
(where $L_b$ should have a factor $1/N$). Much simpler is the case when
$L_b$ fixes the different densities at
the boundaries, whose analysis is carried out sketchily
in the next section simply as an introduction.

\section{Density reservoirs}
\label{sec:z2}

We consider in this section
the Markov process on $\{0,1\}^{[-N,N]}$ with generator $L= L_0
+ L'$, where
denoting by $\eta$ the elements of $\{0,1\}^{[-N,N]}$,\vspace*{1pt}
%
\begin{equation}
\label{z2.1} L_0 f(\eta):=\frac{1}2\sum
_{x=-N}^{N-1} \bigl[f\bigl(\eta^{(x,x+1)}\bigr)-f(
\eta)\bigr]
\end{equation}
with $\eta^{(x,x+1)}(x)= \eta(x+1)$, $\eta^{(x,x+1)} (x+1)=\eta(x)$
and $\eta^{(x,x+1)}(\cdot)=\eta(\cdot)$  elsewhere\vspace*{1pt}
\begin{eqnarray*}
L' f(\eta) &=&\rho_+ \bigl[f\bigl(\eta^{(+,N)}\bigr)-f(\eta)
\bigr] + (1-\rho_+) \bigl[f\bigl(\eta^{(-,N)}\bigr)-f(\eta)\bigr]
\\[2pt]
&&{}+ \rho_- \bigl[f\bigl(\eta^{(+,-N)}\bigr)-f(\eta)\bigr] + (1-\rho_-)
\bigl[f\bigl(\eta ^{(-,-N)}\bigr)-f(\eta)\bigr],
\end{eqnarray*}
where $1\ge\rho_+>\rho_-\ge0$ and $\eta^{+,x}(x)=1$, $\eta
^{+,x}(y)=\eta(y)$, $y\ne x$; analogously,
$\eta^{-,x}(x)=0$, $\eta^{-,x}(y)=\eta(y)$, $y\ne x$.

The process corresponding to $L'$ alone leaves unchanged the
occupations at $|x|<N$ while the equilibrium probabilities
of occupation at $\pm N$ are equal to $\rho_{\pm}$. Since $\rho
_+>\rho_-$, this creates a density gradient and the
full process with generator $L=L_0+L'$ describes the particles flux
determined by the
density gradient. The process is uniformly D\"oblin, in particular,
there is
a unique stationary measure $\mu^{\mathrm{st}}_N$ to which the
process converges
exponentially fast. The averages $\mu^{\mathrm{st}}_N[\eta(x)]$
describe a
linear density profile in agreement with Fick's law. Fluctuations in
the stationary regime are well characterized (\cite{spohn}, and the
large deviations as well, \cite{derrida}).

Denote by $\mu_N$ the initial distribution and by
$\mu_NS_t$ the distribution at time $t$ (i.e.,
the law
at time $t$ of the process with generator $L$ starting from $\mu_N$).
Then, since the process is uniformly D\"oblin, for any positive integer
$N$ there are strictly positive constants $c_N$ and $b_N$ so that\vspace*{1pt}
%
\begin{equation}
\label{z2.2} \bigl \|\mu_NS_t-\mu^{\mathrm{st}}_N
\bigr \| \le c_N \mathrm{e}^{-b_Nt} \qquad\mbox{for any $
\mu_N$ and $t >0$},
\end{equation}
where for any signed measure $\lambda$ on $\{0,1\}^{[-N,N]}$
%
\begin{equation}
\label{z2.3} \|\lambda\| = \sum_{\eta} \bigl |\lambda(
\eta)\bigr |.
\end{equation}
We now prove the following.

%
\begin{them}
\label{thmz2.1}
There are $c$ and $b>0$ independent of $N$ so that for any initial
measure $\mu_N$ and
all $t>0$
%
\begin{equation}\label{z2.4}
\bigl \| \mu_NS_t-\mu^{\mathrm{st}}_N\bigr \| \le c N
\mathrm{e}^{-bN^{-2} t}.
\end{equation}
\end{them}

\begin{pf}
 Let
%
\begin{equation}
\label{z2.5} \mathcal X_N = \bigl\{\und\eta=\bigl(
\eta^{(1)},\eta^{(2)}\bigr)\in\bigl(\{0,1\} \times\{0,1\}\bigr)
^{[-N,N]}\dvtx \eta_{\ne}(x):= \eta^{(1)}(x) - \eta
^{(2)}(x)\ge0, \forall x \bigr\},
\end{equation}
and, for $f\dvtx\mathcal X_N \to\mathbb R$,
\begin{eqnarray*}
\mathcal L_0 f(\und\eta)&:=&\frac{1}2\sum
_{x=-N}^{N-1} \bigl[f\bigl(\und\eta ^{(x,x+1)}
\bigr)-f(\und\eta)\bigr],
\\
\mathcal L' f(\und\eta) &=&\rho_+ \bigl[f\bigl(\und
\eta^{(+,N)}\bigr)- f(\und\eta)\bigr] + (1-\rho_+) \bigl[f\bigl(\und
\eta^{(-,N)}\bigr)-f(\und\eta)\bigr]
\\
&&{}+ \rho_- \bigl[f\bigl(\und\eta^{(+,-N)}\bigr)-f(\und\eta)\bigr] + (1-
\rho_-) \bigl[f\bigl(\und \eta^{(-,-N)}\bigr)-f(\und\eta)\bigr],
\end{eqnarray*}
where $\und\eta^{(+,x)}(x)=(1,1)$, $\und\eta^{(-,x)}(x)=(0,0)$, and
coincide with $\und\eta$ elsewhere, $x=\pm N$.

It is easy to see that $\mathcal L_0$ and $\mathcal L'$ define Markov
generators on $\mathcal X_N$.
Moreover, when acting on functions that depend on only one of the two
entries, $\eta^{(1)}$ or $\eta^{(2)}$,
of $\und\eta$, we see that\vspace*{1.5pt}
$\mathcal L_0+\mathcal L'$ coincide with $L$, and so it defines a
coupling between the processes with generator $L$ starting from two
comparable configurations $\eta^{(1)}$ and $\eta^{(2)}$ ($\eta
^{(1)}(x) \ge\eta^{(2)}(x)$ for all $x$), showing that the $L$-evolution
is attractive in the sense of \cite{li} (i.e., preserves order). In
particular, we may take $\eta^{(1)}\equiv1$ and $\eta^{(2)}\equiv0$
the configurations that are identically 1 and, respectively, 0. Moreover,
$\mathcal L_0$ leaves unchanged
the number of discrepancies which instead may decrease
under the action of~$\mathcal L'$.
Write $\mathbf P$ for the law of the process
starting from $\eta^{(1)}\equiv1$ and $\eta^{(2)}\equiv0$ and call
$\pi(x,t)= \mathbf P[\eta_{\ne}(x,t)=1]$. We then have, recalling
that $\pi(x,0)= 1$ for all $x$,
%
\begin{equation}
\label{z2.6a} \pi(x,t)= 1 - \int_0^t
\bigl(p_s(x,N)\pi(N,t-s)+ p_s(x,-N)\pi (-N,t-s) \bigr)
\,\mathrm{d}s,
\end{equation}
where $p_s(x,y)$ is the probability under the stirring process (with
only one particle) of going from $x$ to $y$ in a time $s$; this is the
same as the probability of a simple random walk whose jumps outside
$[-N,N]$ are suppressed. Indeed, \eqref{z2.6a} follows at once from
the integration by parts formula for the semigroup $\mathcal S_t$
generated by $\mathcal L_0+ \mathcal L'$, with $\mathcal S^0_t$ the
semigroup generated by $\mathcal L_0$, and recalling that the effect of
$\mathcal L'$ is to kill discrepancies at $N$ and $-N$ with rate $1$:
\[
\mathcal S_t (f)= S^0_t (f) +\int
_0^t \mathcal S_{t-s}\bigl(\mathcal
L' S^0_{s}f\bigr) \,\mathrm{d}s,
\]
where $f\dvtx\mathcal X_N \to\mathbb R$.
From \eqref{z2.6a}, we see that
\[
\pi(x,t)= E_x \bigl[ \mathrm{e}^{ - T^*(t)} \bigr],
\]
where $E_x$ is the expectation of the process with transition
probabilities $p_s(x,y)$ and
\[
T^*(t) =\int_0^t ( \mathbf1_{x_s =N}+
\mathbf1_{x_s =-N} ) \,\mathrm{d}s
\]
is the time spent at $\{-N,N\}$ during $[0,t]$. Indeed,
\[
E_x \bigl[ \mathrm{e}^{-T^*(t)} \bigr] =\sum
_{n=0}^\infty\frac{(-1)^n}{n!} \int_0^t
\cdots\int_0^t \mathrm{d}s_1\cdots
\,\mathrm{d}s_n E_x \Biggl[ \prod_{i=1}^n
\{\mathbf1_{x_{s_i} =N}+ \mathbf1_{x_{s_i} =-N}\} \Biggr]
\]
which is the same series which is obtained by iterating \eqref{z2.6a}.

We shall prove that
%
\begin{equation}
\label{z2.6} E_x \bigl[ \mathrm{e}^{ - T^*(t)} \bigr] \le c
\mathrm{e}^{-b N^{-2}t}
\end{equation}
which will then imply
\[
\sum_{x=-N}^N\mathbf P \bigl[
\eta_{\ne}(x,t) =1 \bigr] \le Nc \mathrm {e}^{-b N^{-2}t}
\]
and so \eqref{z2.4}, because $\mu_NS_t$ and $\mu^{\mathrm{st}}_N$ are
squeezed in between the laws of
the marginals of the coupled process.

\textit{Proof of} \eqref{z2.6}.
By an iterative argument, it is enough to show that
\[
\sup_{x\in[-N,N] }E_x \bigl[ \mathrm{e}^{ - \tau}
\bigr] \le p <1,\qquad \tau:= T^*\bigl(N^2\bigr).
\]
But
%
\begin{equation}
\label{z2.7} \inf_{x\in[-N,N]} P_x [ \tau\ge1 ] \ge
\delta>0
\end{equation}
as the probability of reaching $\{-N,N\}$ by time $N^2-1$ is bounded
from below uniformly in the starting
point and the probability of not moving for a unit time interval is
also bounded away from 0. By \eqref{z2.7},
\begin{eqnarray*}
E_x \bigl[ \mathrm{e}^{ - \tau} \bigr] &=& E_x
\bigl[ \mathrm{e}^{ -
\tau}; \tau<1 \bigr] + E_x \bigl[
\mathrm{e}^{ - \tau};\tau\ge1 \bigr]
\\
&\le& 1-P_x[\tau\ge1] + P_x[\tau\ge1]
\mathrm{e}^{-1} \le1 - \delta\bigl(1-\mathrm{e}^{-1}\bigr).
\end{eqnarray*}
\upqed\end{pf}

\section{Main result}
\label{sec:z3}

In this paper, we study the process with generator
$L= L_0+ L_b$, $L_0$ as in \eqref{z2.1}, $L_{b}=L_{b,+}+L_{b,-}$
describes births and deaths near the boundaries. Namely, denoting by
$\eta$ the elements of $\{0,1\}^{[-N,N]}$ and by $f$ functions on $\{
0,1\}^{[-N,N]}$,
%
\begin{eqnarray}
\label{z3.1} L_{b,\pm} f(\eta)&:=& \frac{j}{2N} \sum
_{x\in I_\pm}D_{\pm}\eta(x) \bigl[f\bigl(\eta^{(x)}
\bigr)-f(\eta)\bigr],
\nonumber
\\
D_{+}\eta(x)&=& \bigl(1-\eta(x)\bigr) \eta(x+1)\cdots\eta(N),
\\
D_{-}\eta(x)&=& \eta(x) \bigl(1- \eta(x-1)\bigr)\cdots\bigl(1-
\eta(-N)\bigr),
\nonumber
\end{eqnarray}
where $j>0$ is a parameter of the model, $I_+=\{N-1,N\}$ and $I_-=\{-N,
-N+1\}$ (in \cite{DPTVjsp,DPTVpro,DPTVjsp2} $I_{\pm
}$ consist
of $K$ sites, here we restrict to $K=2$
only for notational simplicity). Thus $L_{b,+}$ adds a particle at rate
$\frac{j}{2N}$ in the last empty site (if any) in $I_+$ while at the
same rate $L_{b,-}$ takes out the first particle (if any) in $I_-$.

Motivations for this model can be found in previous
papers, \cite{DPTVjsp,DPTVpro,DPTVjsp2},
where we have studied the hydrodynamic behavior of the
system and the profile of the stationary measure as $N\to\infty$.
The analysis in the above papers does not say what happens
for the process after
the hydrodynamical regime, that is, at times longer than $ N^2$. This
is the aim of the current paper where we
study the time scale for reaching the stationary regime.

We use the same notation as in the previous section with $S_t=\mathrm
{e}^{Lt}$ and
$\mu_NS_t$, $t\ge0$, the law at time $t$ of the process with
generator $L$ starting from $\mu_N$:
%
\begin{equation}
\label{z3.2} \mu_N S_t [f ] = \mu_N \bigl[
\mathrm{e}^{Lt}f \bigr]=\mu _N \bigl[S_t (f)
\bigr].
\end{equation}

If $j=0$, that is, $L=L_0$ the sets $\{\sum\eta(x)=M\}$, $0\le M \le
2N+1$, are invariant so that
the process is not even ergodic. However, the presence of $L_b$, even
if ``small''
due to the rate $j/2N$, changes drastically the long time behavior of
the system
and it is therefore crucial in the computation of the spectral gap. Our process,
like the one in the previous section, is uniformly D\"oblin; there is
therefore a unique stationary
measure $\mu^{\mathrm{st}}_N$ and \eqref{z2.2} holds in the present
context as well.
We prove the analogue of Theorem~\ref{thmz2.1}.

%
\begin{them}
\label{thmz3.1}
There are $c$ and $b>0$ independent of $N$ so that
\[
\bigl \| \mu_{N}S_t-\mu^{\mathrm{st}}_N\bigr \| \le c N
\mathrm{e}^{-bN^{-2}
t}, \qquad\mbox{for all initial measures $
\mu_{N}$ and all $t >0$.}
\]
\end{them}

Theorem~\ref{thmz3.1} is the main result in this
paper and it will be proved in the next sections.

The rate $N^{-2}$ in the exponent in \eqref{z3.3} cannot be improved,
as can be easily seen
by bounding from below the probability that an initially existing
discrepancy does not
disappear by the time~$N^2$.

The result is in several respects surprising:
the spectral gap in fact scales as $N^{-2}$ just like in the
stirring process (i.e., with $j=0$)
restricted to any of the invariant subspaces $\{\eta\dvtx\sum\eta
(x)= M\}$.
The result says that in a time of the same order the full process
manages to
equilibrate among all the above subsets according to $\mu^{\mathrm{st}}_N$;
also, the
time for this to happen scales in the same way as for the process of the
previous section, where however the birth-death events are not scaled down
with $N$ as in Theorem~\ref{thmz3.1}.\vadjust{\goodbreak}

We do not have sharp information on $\mu^{\mathrm{st}}_N$. In \cite
{DPTVjsp2}, we have proved that
the set $\mathcal M$ of all probability measures on $\{0,1\}^{[-N,N]}$
shrinks after a time of order $N^2$
to a smaller set $\mathcal M_N$ but we have no information on the way
it further shrinks at later times. All measures in $\mathcal M_N$ are
close to a product measure $\ga_N$, meaning that
the expectation of products $\eta(x_1)\cdots\eta(x_n)$ are close
(the accuracy increasing with $N$) to those of $\ga_N$, for all
$n$-tuples of distinct sites $x_i$; $n$ is given, but it can be taken
larger and larger as $N$
increases. We also know that
the expectations $\ga_N[\eta(x)]$ are close to $\rho^{\mathrm
{st}}(x/N)$, where $\rho^{\mathrm{st}}(r)$, $r\in[-1,1]$, is the
stationary solution of the limit hydrodynamic equation; it is an increasing
linear function and $\rho^{\mathrm{st}}(-1)=1-\rho^{\mathrm{st}}(1)>0$.

We thus know that $\mu^{\mathrm{st}}_N$ is close (in the above sense) to
the product measure $\ga_N$, but that is all,
which does not seem detailed enough to apply the usual techniques for
the investigation
of the spectral gap using equilibrium estimates. We proceed
differently, and our proof of Theorem~\ref{thmz3.1} follows along the
lines of the much simpler Theorem~\ref{thmz2.1}. It
relies on a careful analysis of the time evolution, exploiting
stochastic inequalities, as in the previous section.
We thus consider a coupled process on
$\mathcal X_N$ (see \eqref{z2.5}),
which again starts from $\eta^{(1)}(x,0)=1$ and $\eta^{(2)}(x,0)=0$
for all $x\in[-N,N]$. The process is defined in such a
way that the marginal distributions of $\eta^{(1)}$ and
$\eta^{(2)}$ have the law of process with generator $L$. By the
definition of
$\mathcal X_N$,
$\eta^{(1)}\ge\eta^{(2)}$
at all times (order is preserved) and the proof of Theorem~\ref{thmz3.1}
follows from an estimate on the extinction time of
the ``discrepancy configuration''
$\eta_{\ne} = \eta^{(1)} - \eta^{(2)} $. We shall in fact prove
that there are $c$ and $b>0$ independent of $N$ so that
%
\begin{equation}
\label{z3.3} \sum_{x=-N}^N \mathbf P
\bigl[\eta_{\ne}(x,t)=1\bigr] \le c N \mathrm {e}^{-bN^{-2} t}.
\end{equation}

\section{The coupled process}
\label{sec:z4}

Throughout the sequel, we shall use the following.

\begin{Notation*}
$\eps:=N^{-1}$; for $\und\eta= (\eta^{(1)},
\eta^{(2)}) \in\mathcal X_N$ as defined in \eqref{z2.5}, and
$x \in[-N,N]$,
%
\begin{eqnarray}
\label{z4.1} \eta_{\ne}(x)&=& \eta^{(1)}(x)-
\eta^{(2)}(x),
\nonumber
\\
\eta_1(x)&=& \eta ^{(1)}(x) \eta^{(2)}(x),
\\
\eta_0(x)&=&\bigl(1-\eta^{(1)}(x)\bigr) \bigl(1-
\eta^{(2)}(x)\bigr),
\nonumber
\end{eqnarray}
$\eta_{\ne},\eta_1,\eta_0$ are all in $\{0,1\}^{[-N,N]}$ and $\eta
_{\ne}+\eta_1+\eta_0\equiv1$.
Thus, \eqref{z4.1} establishes a one-to-one correspondence between
$\mathcal X_N$ and $\{\ne,1,0\}^{[-N,N]}$.
By an abuse of notation, we shall denote again by $\und\eta$ the
elements of $\{\ne,1,0\}^{[-N,N]}$, thinking of $\eta_{\ne},\eta
_1,\eta_0$ as functions of $\und\eta$. We may then say that a $\ne
$, $1$ or $0$-particle is
at $x$ according to the value of $\und\eta(x)$.
\end{Notation*}

\begin{Definition*}
Call $L'_0$ the stirring generator acting on functions on $\mathcal
X_N$ (defined as in \eqref{z2.1} with $\eta$ replaced by $\und\eta
$) and let
$L_c= L'_0+\frac{j}{2N}L_1$, $L_1=L_r+L_l$, be the generator acting on
functions on $\mathcal X_N$, where
$L_rf$ is defined as
%
\begin{eqnarray}
\label{z4.2} L_rf(\und\eta)&=&\sum_{i=N-1}^N
D(\und\eta,i) \bigl[f\bigl(\und\eta^{\ne
,1,i}\bigr) -f(\und\eta)\bigr]
\nonumber
\\[-8pt]
\\[-8pt]
&&{} + A(\und\eta,N) \bigl[f\bigl(\und\eta^{\ne,1,N;0,\ne,N-1}\bigr)-f(\und\eta)
\bigr] +\sum_{i=N-1}^N B(\und\eta,i)\bigl[f
\bigl(\und\eta^{0,1,i}\bigr)-f(\und\eta)\bigr]
\nonumber
\end{eqnarray}
and where $\und\eta^{a,b,i}$ changes from $a$ to $b$ the value of
$\und\eta$ at site $i$ if $\und\eta(i)=a$, and
$\und\eta^{a,b,i}=\und\eta$ otherwise, and $\und\eta^{\ne
,1,N;0,\ne,N-1}=(\und\eta^{\ne,1,N})^{0,\ne,N-1}$,
\begin{eqnarray*}
\nonumber
D(\und\eta,N) &=& \eta_{\ne}(N)\bigl[1-\eta_0(N-1)
\bigr],\qquad D(\eta,N-1)=\eta _{\ne}(N-1)\eta_1(N),
\\
A(\und\eta,N) &=& \eta_{\ne}(N) \eta_0(N-1),
\\
B(\und\eta,N) &=& \eta_0(N),\qquad B(\eta,N-1)=
\eta_0(N-1)\eta_1(N).
\end{eqnarray*}
\end{Definition*}

Thus, $L_r$ describes three types of events all occurring in $I_+$:
\begin{itemize}
\item$D$-events: a $\ne$-particle becomes a $1$-particle.

\item$A$ events: a $\ne$-particle becomes a $1$-particle and simultaneously
a $0$-particle becomes a $\ne$-particle.

\item$B$-events: a $0$-particle becomes a $1$-particle.
\end{itemize}

$L_l$ is defined analogously by changing $I_+$ into $I_-$ and $\eta_0$
with $\eta_1$.
One can easily check that
%
\begin{equation}
\label{z4.3} L_c f = Lg, \qquad\mbox{whenever $f(\und\eta) = g
\bigl(\eta^{(i)}\bigr)$, $i=1,2$},
\end{equation}
$L$ the generator in Section~\ref{sec:z3}. Thus, the process generated
by $L_c$ is a coupling of two
processes both with generator $L$ and   $L$ preserves order (this is
just the standard basic coupling, as in \cite{li}; see also
Proposition~3.1 of \cite{DPTVjsp2}).

\section{Graphical construction}
\label{sec:z5}

Following the so-called Harris graphical construction, we realize the
coupled process in a
probability space $(\Omega, \mathcal{F}, P)$ where several independent Poisson
processes are defined.

\begin{Definition*} The probability space $(\Omega, \mathcal{F}, P)$. The elements
$\omega\in\Om$ have the form
\[
\omega= \bigl( \und t^{(x)}, x\in[-N,N-1]; \und t^{(A, \pm N)} \und
t^{(D, \pm N)};\und t^{(D,\pm(N-1))}; \und t^{(B, \pm N)};\und
t^{(B,\pm(N-1))} \bigr),
\]
where each entry is a sequence in $\mathbb R_+$ whose elements are
interpreted as times. Under $P$, the
entries are independent Poisson processes: each one of the $\und
t^{(x)}$ has intensity $1/2$,
and all the others have each intensity $\eps j/2$.
\end{Definition*}

With probability 1, all times are different from each other
and there are finitely many
events in a compact. For any\vadjust{\goodbreak} such $\omega\in\Om$,
we construct piecewise constant functions $\eta_1(x,t;\omega)$, $\eta
_0(x,t;\omega)$, $\eta_{\ne}(x,t;\omega)$,
as follows.\vspace*{-1pt}
The jump times are a subset of the events in the above Poisson
processes, more specifically
at the times $t=t^{(x)}_{n}$ we exchange
the content of the sites $x$ and $x+1$ (i.e., we do a
stirring at $(x,x+1)$); the other jumps are:\vspace*{-2pt}
\begin{itemize}
\item
At the times $t=t^{(A,\pm N)}_{n}$, the configuration is updated only
if $\eta_{\ne}(\pm N,t^-)=1$,  $\eta_{0}(\pm(N-1),t^-)=1$ and the
new configuration has $\eta_{\ne}(\pm(N-1),t^+)=1$ and $\eta
_{1}(\pm N,\allowbreak  t^+)=1$; the values at other sites remain unchanged.

\item At the times $t=t^{(D,\pm N)}_{n}$, the configuration is updated only
if $\eta_{\ne}(\pm N,t^-)=1$ and $\eta_{0}(\pm(N-1),t^-)=0$, the
new configuration has
$\eta_{1}(\pm N,t^+)=1$; the values at other sites unchanged.

\item At the times $t=t^{(D,\pm(N-1))}_{n}$, the configuration is
updated only
if $\eta_{\ne}(\pm(N-1),t^-)=1$ and $\eta_{1}(\pm N,t^-)=1$; the
new configuration has
$\eta_{1}(\pm(N-1),t^+)=1$; the values at other sites unchanged.

\item
At the times $t=t^{(B,\pm N)}_{n}$, the configuration is updated only
if $\eta_{0}(\pm N,t^-)=1$; the new configuration has
$\eta_{1}(\pm N,t^+)=1$; the values at other sites unchanged.

\item
At the times $t=t^{(B,\pm(N-1))}_{n}$, the configuration is updated only
if $\eta_{1}(N,t^-)=1$ and $\eta_{0}(\pm(N-1),t^-)=1$; the new
configuration has
$\eta_{1}(\pm(N-1),t^+)=1$; the values at other sites unchanged.\vspace*{-1pt}
\end{itemize}

We take initially $\eta_{\ne}(x,0)=1$ for all $x$,
then the variables $\und\eta(x,t;\omega)$ defined as above on $(\Om,P)$
have the law of the coupled process defined in Section~\ref{sec:z4}.\vspace*{-2pt}

\begin{Definition*}
Labeling the discrepancies. By realizing the
process in
the space $(\Omega, \mathcal{F}, P)$, we can actually follow the discrepancies in time.
Indeed consider the discrepancy
initially at a site $z \in[-N,N]$. Then the discrepancy will move
following the marks of $\omega$. Namely, it moves
at the stirring times, that is, it jumps from $x$ to $x+1$
(or from $x+1$ to $x$) at the times $t\in\und t^{(x)}$. Moreover, it
jumps from $N$ to $N-1$ at the times in
$\und t^{(A,N)}$ (if $\eta_0(N-1)=1$) and analogously from $-N$ to
$-N+1$ at the times in
$\und t^{(A,-N)}$ (if $\eta_1(-N+1)=1$). Finally, we say that the
discrepancy dies (and goes to the state $\emptyset$)
at the times $\und t^{(D,\pm N)}$, $\und t^{(D,\pm(N-1))}$ (if the
conditions for the event are satisfied, as explained in the previous
paragraphs).\vspace*{-2pt}
\end{Definition*}

We thus label the initial discrepancies by assigning with uniform
probability a
label in $\{1,\ldots,2N+1\}$ to each site in $[-N,N]$ and call $(z_1,\ldots
,z_{2N+1})$
the sites corresponding to the labels $1,\ldots,2N+1$. This is done
independently
of $\omega$ and by an abuse of notation we still denote by $P$ the
joint law of $\omega$ and the labeling. Since initially all sites
are occupied by discrepancies, we may
interpret $z_i$ as the position at time $0$
of the discrepancy with label $i$. In particular at time $0$, the
probability that
$z_i=x$ is equal to $1/(2N+1)$.
Given $\omega\in\Om$, we
follow the motion of the
labeled discrepancies as described above and define accordingly
the variables $z_i(t,\omega)$ which take values in $\{[-N,N]\cup
\emptyset\}$. Thus, the set
$Z(t,\omega)$ of all $z_i(t,\omega)\ne\emptyset$ is equal to $\{x\dvtx
\eta_{\ne}(x,t;\omega) =1\}$,
so that\vspace*{-2pt}
%
\begin{eqnarray}
\label{z5.1} P \biggl[ \sum_x
\eta_{\ne}(x,t)>0 \biggr] &=& P \bigl[ \mbox{there is $i\dvtx
z_i(t,\omega) \ne\emptyset$} \bigr]
 \le \sum_i P\bigl[z_i(t,\omega)
\ne\emptyset\bigr]
\nonumber
\\[-9pt]
\\[-9pt]
&= & (2N+1) P\bigl[z_1(t,\omega) \ne\emptyset\bigr],
\nonumber
\end{eqnarray}
the last equality by symmetry.\vadjust{\goodbreak}

Obviously, $P[z_1(t,\omega) \ne\emptyset]$ does not depend on the
labels of the other $z$-particles
so that we may and shall describe the system in terms of a random walk
$z_t=z_1(t,\omega)$ in a random
environment $\eta_t\in\{\ne,0,1\}^{[-N,N]\setminus z_t}$ when
$z_t\ne\emptyset$ (i.e., it is alive); when $z_t= \emptyset$ then
$\eta_t\in\{\ne,\allowbreak  0,1\}^{[-N,N]}$, but since we want to study
$P[z_1(t,\omega) \ne\emptyset]$ what happens after the death of $z$
is not relevant.

We have reduced the problem to the analysis of
the extinction time
of a random walk in a random environment: the problem looks now
very similar to the one considered in \cite{dptvrw}, the only
difference being that the environment
has a more complex structure with three rather than two states per
site. But the procedure is essentially the same
as we briefly sketch in the sequel.

\section{The auxiliary random walk process}
\label{sec:6}

Once the initial condition $(z,\eta^*)$ has been fixed, we can
consider an auxiliary time dependent Markov process $(\tilde z_t)$ as
in \cite{dptvrw}, whose extinction time has the same law as that of
the true process $(z_1(t))$ of the previous section. The transition
rates for
$\tilde z_t$ are given by the conditional expectation of the transition
rates of $(z_1(\cdot))$ conditioned on $z_1(t)$. Thus, they depend
on the law of the full process, and hence on the initial datum $(z,\eta^*)$.
This time dependent generator $\mathcal L_t$ is given in \eqref{6.3}
below, and satisfies
\[
\tilde E_{z}\bigl[\mathcal L_t f(\tilde z_t)
\bigr]=E_{z,\eta^*} \bigl[L \phi \bigl(z_1(t),\eta_t
\bigr) \bigr]=\frac{\mathrm{d}}{\mathrm{d}t} E_{z,\eta^*} \bigl[ \phi \bigl(z_1(t),
\eta_t\bigr) \bigr],
\]
where $\phi(z,\eta)=f(z)$ and $f\dvtx \La_N\cup\emptyset\to\mathbb R$.

Since
\begin{eqnarray*}
L_{r} \phi&=& \frac{j}{2N} \bigl\{\mathbf1_{z=N}
\bigl(1- \eta_0(N-1)\bigr) \bigl[f(\emptyset)-f(N)\bigr]
\\
&&\qquad{} +\mathbf1_{z=N-1} \eta_1(N) \bigl[f(\emptyset)-f(N-1)\bigr]
\bigr\}
\\
&&{}+ \frac{j}{2N} \mathbf1_{z=N} \eta_0(N-1)
\bigl[f(N-1)-f(N)\bigr],
\\
L_l \phi&=&\frac{j}{2N} \bigl\{\mathbf1_{z=-N}
\bigl(1-\eta_1(-N+1)\bigr) \bigl[f(\emptyset)-f(-N)\bigr]
\\
&&\qquad{}+
\mathbf1_{z=-N+1} \eta_0(-N) \bigl[f(\emptyset)-f(-N+1)\bigr]
\bigr\}
\\
&&{}+ \frac{j}{2N} \mathbf1_{z=-N} \eta_1(-N+1)
\bigl[f(-N+1)-f(-N)\bigr],
\end{eqnarray*}
we set
%
\begin{eqnarray}
d(N,t) &=& \frac{j}{2N} E_{z_0,\eta^*}\bigl[1-\eta_0(N-1,t) |
z_t=N\bigr],
\nonumber
\\
d(N-1,t) &=& \frac{j}{2N} E_{z_0,\eta^*}\bigl[\eta_1(N,t) |
z_t=N-1\bigr],
\nonumber
\\
d(-N,t) &=& \frac{j}{2N} E_{z_0,\eta^*}\bigl[\bigl(1-
\eta_1(-N+1,t)\bigr) | z_t=-N\bigr],
\nonumber
\\
\label{6.1}d(-N+1,t) &=& \frac{j}{2N} E_{z_0,\eta_0}\bigl[
\eta_0(-N,t) | z_t=-N+1\bigr],
\\
a(N,t) &=& \frac{j}{2N} E_{z_0,\eta^*}\bigl[ \eta_0(N-1,t) |
z_t=N\bigr],
\nonumber
\\
\label{6.2}a(-N,t) &=& \frac{j}{2N} E_{z_0,\eta^*}\bigl[
\eta_1(-N+1,t) | z_t=-N\bigr],
\end{eqnarray}
and $d(z,t)=0$ if $|z|<N-1$. Thus, for $t\ge0$, we have
%
\begin{eqnarray}
\label{6.3} \mathcal L_t f(z) &=& \mathcal L^0 f (z)+
d(z,t) \bigl[f(\emptyset)-f(z)\bigr] + \mathbf1_{z=N}a (N,t)
\bigl[f(N-1)-f(N)\bigr]
\nonumber
\\[-8pt]
\\[-8pt]
&&{}+ \mathbf1_{z=-N}a (-N,t) \bigl[f(-N+1)-f(-N)\bigr].
\nonumber
\end{eqnarray}
The process $\tilde z_t$ is a simple random walk with extra jumps from
$N$ to $N-1$ and $-N$ to $-N+1$
with time-dependent intensity $a(\pm N,t)$; moreover, it has death rate
$d(z,t)$ (rate to go to $\emptyset$). Observe that
\[
d(z,t) \ge\frac{j}{2N} E_{z_0,\eta^*}\bigl[\eta_1(N-1,t) |
z_t=N\bigr]\mathbf1_{z=N},
\]
and the analysis becomes very similar to the case treated in \cite{dptvrw}.
From the same argument leading to Theorem~1 therein, we have that for
any initial configuration $\eta^*$ and $z_0$:
\[
P\bigl[z_1(t)\neq\emptyset\bigr]\le c \mathrm{e}^{-b N^{-2}t},
\]
which completes the proof.


\section*{Acknowledgements}

The research has been partially supported by PRIN 2009 (2009TA2595-002).
The research of D.~Tsagkarogiannis is partially supported by the FP7-REGPOT-2009-1 project
``Archimedes Center for Modeling, Analysis and Computation'' (under
grant agreement no. 245749).
M. E. Vares is partially supported by CNPq grants PQ 304217/2011-5 and 474233/2012-0.


%

\printhistory

\begin{thebibliography}{11}

\bibitem{bertini}
%
\begin{barticle}[mr]
\bauthor{\bsnm{Bertini},~\bfnm{L.}\binits{L.}},
\bauthor{\bsnm{De Sole},~\bfnm{A.}\binits{A.}},
\bauthor{\bsnm{Gabrielli},~\bfnm{D.}\binits{D.}},
\bauthor{\bsnm{Jona-Lasinio},~\bfnm{G.}\binits{G.}} \AND
\bauthor{\bsnm{Landim},~\bfnm{C.}\binits{C.}}
(\byear{2006}).
\btitle{Non equilibrium current fluctuations in stochastic lattice gases}.
\bjournal{J. Stat. Phys.}
\bvolume{123}
\bpages{237--276}.
\bid{doi={10.1007/s10955-006-9056-4}, issn={0022-4715}, mr={2227084}}
\end{barticle}
%
\bptok{imsref}%
\endbibitem

\bibitem{bodineau}
%
\begin{barticle}[mr]
\bauthor{\bsnm{Bodineau},~\bfnm{T.}\binits{T.}} \AND
\bauthor{\bsnm{Derrida},~\bfnm{B.}\binits{B.}}
(\byear{2006}).
\btitle{Current large deviations for asymmetric exclusion processes
with open boundaries}.
\bjournal{J. Stat. Phys.}
\bvolume{123}
\bpages{277--300}.
\bid{doi={10.1007/s10955-006-9048-4}, issn={0022-4715}, mr={2227085}}
\end{barticle}
%
\bptok{imsref}%
\endbibitem

\bibitem{bdl}
%
\begin{barticle}[mr]
\bauthor{\bsnm{Bodineau},~\bfnm{T.}\binits{T.}},
\bauthor{\bsnm{Derrida},~\bfnm{B.}\binits{B.}} \AND
\bauthor{\bsnm{Lebowitz},~\bfnm{J.~L.}\binits{J.L.}}
(\byear{2010}).
\btitle{A diffusive system driven by a battery or by a smoothly
varying field}.
\bjournal{J. Stat. Phys.}
\bvolume{140}
\bpages{648--675}.
\bid{doi={10.1007/s10955-010-0012-y}, issn={0022-4715}, mr={2670735}}
\end{barticle}
%
\bptok{imsref}%
\endbibitem


\bibitem{DPTVjsp}
%
\begin{barticle}[mr]
\bauthor{\bsnm{De Masi},~\bfnm{A.}\binits{A.}},
\bauthor{\bsnm{Presutti},~\bfnm{E.}\binits{E.}},
\bauthor{\bsnm{Tsagkarogiannis},~\bfnm{D.}\binits{D.}} \AND
\bauthor{\bsnm{Vares},~\bfnm{M.~E.}\binits{M.E.}}
(\byear{2011}).
\btitle{Current reservoirs in the simple exclusion process}.
\bjournal{J. Stat. Phys.}
\bvolume{144}
\bpages{1151--1170}.
\bid{doi={10.1007/s10955-011-0326-4}, issn={0022-4715}, mr={2841919}}
\end{barticle}
%
\bptok{imsref}%
\endbibitem

\bibitem{DPTVpro}
%
\begin{barticle}[mr]
\bauthor{\bsnm{De Masi},~\bfnm{Anna}\binits{A.}},
\bauthor{\bsnm{Presutti},~\bfnm{Errico}\binits{E.}},
\bauthor{\bsnm{Tsagkarogiannis},~\bfnm{Dimitrios}\binits{D.}} \AND
\bauthor{\bsnm{Vares},~\bfnm{Maria~E.}\binits{M.E.}}
(\byear{2012}).
\btitle{Truncated correlations in the stirring process with births and deaths}.
\bjournal{Electron. J. Probab.}
\bvolume{17}
\bpages{no. 6, 35~pp}.
\bid{doi={10.1214/EJP.v17-1734}, issn={1083-6489}, mr={2878785}}
\end{barticle}
%
\bptok{imsref}%
\endbibitem

\bibitem{DPTVjsp2}
%
\begin{barticle}[mr]
\bauthor{\bsnm{De Masi},~\bfnm{Anna}\binits{A.}},
\bauthor{\bsnm{Presutti},~\bfnm{Errico}\binits{E.}},
\bauthor{\bsnm{Tsagkarogiannis},~\bfnm{Dimitrios}\binits{D.}} \AND
\bauthor{\bsnm{Vares},~\bfnm{Maria~Eulalia}\binits{M.E.}}
(\byear{2012}).
\btitle{Non-equilibrium stationary states in the symmetric simple
exclusion with births and deaths}.
\bjournal{J. Stat. Phys.}
\bvolume{147}
\bpages{519--528}.
\bid{doi={10.1007/s10955-012-0481-2}, issn={0022-4715}, mr={2923327}}
\end{barticle}
%
\bptok{imsref}%
\endbibitem

\bibitem{dptvrw}
%
\begin{bmisc}[auto:STB|2014/05/28|10:36:42]
\bauthor{\bsnm{De Masi},~\bfnm{A.}\binits{A.}},
\bauthor{\bsnm{Presutti},~\bfnm{E.}\binits{E.}},
\bauthor{\bsnm{Tsagkarogiannis},~\bfnm{D.}\binits{D.}} \AND
\bauthor{\bsnm{Vares},~\bfnm{M.~E.}\binits{M.E.}}
(\byear{2014}).
\bhowpublished{Extinction time for a random walk in a random
environment. \textit{Bernoulli}. To appear}.
\end{bmisc}
%
\bptok{imsref}%
\endbibitem



\bibitem{derrida}
%
\begin{barticle}[mr]
\bauthor{\bsnm{Derrida},~\bfnm{B.}\binits{B.}},
\bauthor{\bsnm{Lebowitz},~\bfnm{J.~L.}\binits{J.L.}} \AND
\bauthor{\bsnm{Speer},~\bfnm{E.~R.}\binits{E.R.}}
(\byear{2002}).
\btitle{Large deviation of the density profile in the steady state of
the open symmetric simple exclusion process}.
\bjournal{J. Stat. Phys.}
\bvolume{107}
\bpages{599--634}.
\bid{doi={10.1023/A:1014555927320}, issn={0022-4715}, mr={1898851}}
\end{barticle}
%
\bptok{imsref}%
\endbibitem
\bibitem{li}
%
\begin{bbook}[mr]
\bauthor{\bsnm{Liggett},~\bfnm{Thomas~M.}\binits{T.M.}}
(\byear{1985}).
\btitle{Interacting Particle Systems}.
\bseries{Grundlehren der Mathematischen Wissenschaften [Fundamental
Principles of Mathematical Sciences]}
\bvolume{276}.
\blocation{New York}:
\bpublisher{Springer}.
\bid{doi={10.1007/978-1-4613-8542-4}, mr={0776231}}
\end{bbook}
%
\bptok{imsref}%
\endbibitem

\bibitem{lu}
%
\begin{barticle}[mr]
\bauthor{\bsnm{Lu},~\bfnm{Sheng~Lin}\binits{S.L.}} \AND
\bauthor{\bsnm{Yau},~\bfnm{Horng-Tzer}\binits{H.-T.}}
(\byear{1993}).
\btitle{Spectral gap and logarithmic {S}obolev inequality for
{K}awasaki and {G}lauber dynamics}.
\bjournal{Comm. Math. Phys.}
\bvolume{156}
\bpages{399--433}.
\bid{issn={0010-3616}, mr={1233852}}
\end{barticle}
%
\bptok{imsref}%
\endbibitem

\bibitem{spohn}
%
\begin{barticle}[mr]
\bauthor{\bsnm{Spohn},~\bfnm{Herbert}\binits{H.}}
(\byear{1983}).
\btitle{Long range correlations for stochastic lattice gases in a
nonequilibrium steady state}.
\bjournal{J. Phys. A}
\bvolume{16}
\bpages{4275--4291}.
\bid{issn={0305-4470}, mr={0732737}}
\end{barticle}
%
\bptok{imsref}%
\endbibitem
\end{thebibliography}
\end{document}